\documentclass[12pt,leqno]{amsart}  
\usepackage{amsmath,amstext,amsthm,amssymb}
\usepackage[colorlinks,citecolor=red,pagebackref,hypertexnames=false]{hyperref}
\usepackage[backrefs]{amsrefs}

\def\norm#1.#2.{\lVert#1\rVert_{#2}}
\def\Norm#1.#2.{\bigl\lVert#1\bigr\rVert_{#2}}
\def\NOrm#1.#2.{\Bigl\lVert#1\Bigr\rVert_{#2}}
\def\NORm#1.#2.{\biggl\lVert#1\biggr\rVert_{#2}}
\def\NORM#1.#2.{\Bigl\lVert#1\Bigr\rVert_{#2}}

\def\ip#1,#2,{\langle #1,#2\rangle}
\def\Ip#1,#2,{\bigl\langle#1,#2\bigr\rangle}
\def\IP#1,#2,{\Bigl\langle#1,#2\Bigr\rangle}

\def\mid{\,:\,}

\def\abs#1{\lvert#1\rvert}
\def\Abs#1{\bigl\lvert#1\bigr\rvert}

\numberwithin{equation}{section}

\newtheorem{proposition}[equation]{Proposition} 
\newtheorem{theorem}[equation]{Theorem}

\theoremstyle{remark}

\newtheorem*{Acknowledgment}{Acknowledgment}

\def\eqdef{\stackrel{\mathrm{def}}{{}={}}}

\title[Commutators and Reisz Potentials] {Commutators with Reisz Potentials \\ in One and Several Parameters  }

\author[M.~T.~Lacey]{Michael T. Lacey\\ Georgia Institute of Technology}

\address{Michael Lacey\\
School of Mathematics\\
Georgia Institute of Technology\\
Atlanta,  GA 30332 USA}

\email{lacey@math.gatech.edu}

\thanks{Research supported in part by an NSF grant and a fellowship from the Guggenheim 
Foundation. The hospitality of the University of British Columbia is gratefully acknowledged.}


\begin{document} \parskip=11pt

\begin{abstract} 
Let $\operatorname M_b$ be the operator 
of pointwise multiplication by $b$, that is $\operatorname M_b f=bf$. Set $[\operatorname A,\operatorname B]={} 
\operatorname A \operatorname B-\operatorname B \operatorname A$.  
The Reisz potentials are the operators 
\begin{equation*}  
\operatorname R_\alpha f(x)=\int f(x-y)\frac{dy}{\abs y ^{\alpha} },\qquad 0<\alpha<1.
\end{equation*}
They map $L^p\mapsto L^q$, for $1-\alpha+\frac1q=\frac1p$, a fact we shall take 
for granted in this paper.  A Theorem 
of Chanillo \cite{MR84j:42027}   states that one has the equivalence 
\begin{equation*}  
\norm [\operatorname M_b,\operatorname R_\alpha].p\to q.\simeq \norm b.\operatorname{BMO}.
\end{equation*}
with the later norm being that of the space of functions of bounded mean oscillation.  
We discuss a proof of this result in a discrete setting, and extend part of the equivalence 
above to the higher parameter setting. 
\end{abstract}

\maketitle  

\section{Introduction: One Parameter } \label{s.reisz} 

Our topic is norm bounds on commutators of different operators with the operation of multiplication by a function. 
Chanillo \cite{MR84j:42027} proved that commutators with Reisz potentials characterize the function space $\operatorname{BMO}$.  We 
are concerned with a proof of this result, and extensions  to multi-parameter situations.  

To set notation, let $\operatorname H$ be the Hilbert transform, that is 
\begin{equation*}
\operatorname H f(x)=\text{p.v.}\int f(x-y)\frac{dy}y. 
\end{equation*}
Let $\operatorname M_b$ be the operator 
of pointwise multiplication by $b$, that is $\operatorname M_b f=bf$ and $[\operatorname A,\operatorname B]={} 
\operatorname A \operatorname B-\operatorname B \operatorname A$.  A classical result
states that 
\begin{equation}  \label{e.nehari}
\norm [\operatorname M_b, \operatorname H].p\to p.\simeq \norm b.\operatorname{BMO}.,\qquad 1<p<\infty. 
\end{equation}
The latter space is that of functions of bounded mean oscillation, the dual to the real valued 
Hardy space $H^1$. The dyadic version of this space is defined in (\ref{e.BMO}).

The history of this result goes back to the characterization of the boundedness of Hankel operators due 
to Nehari \cite{nehari}.  In the purely harmonic analysis setting, Coifmann, Rochberg and Weiss \cite{MR54:843} provided an extensive 
study of this equivalence. 

Set the Reisz potentials to be 
\begin{equation}  \label{e.reisz}
\operatorname R_\alpha f(x)=\int f(x-y)\frac{dy}{\abs y ^{\alpha} },\qquad 0<\alpha<1.
\end{equation}
These operators map $L^p\mapsto L^q$, for $1-\alpha+\frac1q=\frac1p$, a fact we shall take 
for granted in this paper. We are interested in the Theorem 
of Chanillo \cite{MR84j:42027}.

\begin{theorem} \label{t.chanillo}
For $1-\alpha+\frac1q=\frac1p$, and $1<p<q<\infty$ we have 
\begin{equation}  \label{e.chanillo}
\norm [\operatorname M_b,\operatorname R_\alpha].p\to q.\simeq \norm b.\operatorname{BMO}.
\end{equation}
\end{theorem}

The method of proof introduced by Chanillo  \cite{MR84j:42027} is to dominate the sharp
function of the commutator, a method that has been 
extended by a variety of authors in different settings, see \cites{MR2004a:42010,
MR2002g:42016
,
MR80j:42034
,
MR96j:42006
,
MR96k:47089 }.  
We give a new proof,  showing that the  commutator with the Reisz potential is 
a sum  of paraproducts.  See the  next section for a definition of paraproducts.

The rationale for this new proof is an extension of part Chanillo's result to a higher parameter setting, 
motivated in part by an extension of the Nehari theorem to higher parameter settings in papers of 
Ferguson and Lacey \cite{MR1961195} and Lacey and Terwelleger \cite{witherin}.  Also see the recent papers 
of Muscalu, Pipher, Tao and Thiele \cites{camil,camil2}. 

Theorem~\ref{t.reiszmany} is the main result of this paper.  It is a partial extension of the one parameter 
result above, in that the upper bound on the commutator is established.  This upper bound is in terms of the 
product $\operatorname {BMO} $ norm of the symbol $b$.  Product $\operatorname {BMO}$ is the one 
identified by S.-Y.~Chang and R.~Fefferman \cites{MR82a:32009,MR86g:42038}, a definition of 
which we will recall below. The lower bound on the 
commutator norm may not be true.  See Section~\ref{s.lower}.

 We restrict attention to discrete forms of the Reisz potentials.   In situations such as this one, 
 it permits one to concentrate on the most essential parts of the proofs, see e.g.~\cite{MR1934198} 
 which is just one reference that is closely associated with the themes of this paper.

 Appropriate averaging procedures will 
 permit one to recover the continuous analogs, but we omit this argument, as it is well represented in the 
 literature. For different versions of this argument, see \cites{MR1783613,MR1964822,MR2000m:42016}.
 In addition, the Reisz potentials on higher dimensional spaces are also frequently considered. 
 The author is not aware of any reason why the methods of this paper will not extend to this 
 level of generality;  it is not pursued as it would complicate our presentation.
 
 \begin{Acknowledgment} 
 We thank the referee for a careful reading of the manuscript.
 \end{Acknowledgment}
 
\section{The One Parameter Statement} \label{s.one}

\subsection{Haar Functions and Paraproducts}

 The dyadic intervals are 
 \begin{equation*}
 \mathcal D {}\eqdef{}\{ [j2^k,(j+1)2^k)\mid j,k\in \mathbb Z\}. 
 \end{equation*}
 Each dyadic interval $I  $ is a union of 
 its left and right halves $I _{-} $, and $I_+ $ respectively.  The Haar function $h_I $ adapted to $I $ is 
 \begin{equation}  \label{e.Haar}
 h_I {}\eqdef{}\abs I ^{-1/2} ( -\mathbf 1 _{I_-}+\mathbf 1 _{I+} ).  
 \end{equation}
 We will also denote the Haar functions as $h^{0}_I $, setting 
 \begin{equation}  \label{e.h1}
 h^1 _{I} = \abs I ^{-1/2} \mathbf 1 _{I}.  
 \end{equation} 
 Thus, $h^{0}_I $ has integral zero, while $h^1_I $ is a multiple of an indicator function. 
 
 It is an essential fact that the Haar functions form an unconditional basis for $L^p$, in particular 
 \begin{equation}  \label{e.littlewoodpaley}
 \norm f.p.\simeq{} \NOrm \Bigl[ \sum _{I\in \mathcal D}  \abs{\ip f,h_I, h^1_I }^2
 		\Bigr] ^{1/2} .p. 
 \end{equation}
  
 Define the dyadic $\operatorname{BMO}$ 
 semi norm by 
 \begin{equation}  \label{e.BMO}
 \norm f.\operatorname{BMO}. {}\eqdef{}\sup _{J\in\mathcal D} \Bigl[ \frac1{\abs J }\sum _{I\subset J} \abs{\ip f,h_I, }^2\Bigr] ^{1/2}.
 \end{equation}

 The Haar paraproducts are 
 \begin{equation}  \label{e.Haar-para-one}
 \operatorname B (f _{1},f _{2})\eqdef{}\sum _{I\in\mathcal D}
 \frac {\ip f _{1},h _{I},}{\sqrt{\abs{I} }}  \ip f _{2},h^1 _{I}, \, h_I
 \end{equation}
 It is critical that there is exactly one function which is a multiple of the identity. 
 We take as a given the fundmental fact about the boundedness of these operators.

\begin{theorem} \label{t.Haar-one}
We have 
 \begin{equation}  \label{e.Haar-one}
\norm  \operatorname B (f_1,\cdot).p. {}\simeq{} 
\norm f_1.\operatorname{BMO}   . ,\qquad 1<p<\infty.
 \end{equation}
\end{theorem} 

This theorem goes back to the work of Coifman and Meyer \cites{MR518170,MR511821,MR0380244}. 
It plays a critical role in the $\operatorname T1$ theorem of David and Journ\'e \cite{MR763911}. 
See for instance the discussion in the text of E.M.~Stein \cite{MR1232192}.  An analogous result 
in the higher parameter situation will be stated and proved  in the next section.

 We shall also appeal to some operators, related to, but not as central,  the paraproducts. 
 Define 
 \begin{equation}  \label{e.operatorC}
 \operatorname C (f _{1},f _{2})\eqdef{}\sum _{I\in\mathcal D}
 \abs{I} ^{-1/2}h_I\prod _{j=1}^{2} \ip f _{j},h _{I},
 \end{equation}
 Notice that every Haar function that appears has zero integral.  Therefore, we can estimate 
 \begin{align*}
 \norm \operatorname C (f _{1},f _{2}).p. & {}\simeq{} 
 	\NOrm \Bigl[ \sum _{I\in \mathcal D}  \abs {\abs{I} ^{-1/2} \prod _{j=1}^{2} \ip f _{j}, h_I , h^1_I }^2 
	\Bigr]  ^{1/2} .p.
\\& {}\le{}  
	\sup _{I\in \mathcal D} \frac{ \abs { \ip f_1,h_I ,} } { \sqrt{\abs{I}} }  
	\NOrm \Bigl[ \sum _{I\in \mathcal D} \abs{\ip f_2,h_I, h^1_I }^2
 		\Bigr] ^{1/2} .p. 
\\ & {}\lesssim{} \sup _{I\in  \mathcal D} \frac{ \abs { \ip f_1,h_I ,} } { \sqrt{\abs{I}} }  \norm f_2.p. 
 \end{align*}
 We will not have recourse to these operators until the next section.

 Define Haar projections to a particular scale by 
 \begin{equation}  \label{e.haarprojections}
 \operatorname P _n f {}\eqdef{}\sum _{\abs I=2^n} \ip f,h_I,h_I . 
 \end{equation}
 And define a related paraproduct by 
 \begin{equation}  \label{e.operatorD}
 \operatorname D_k (f_1,f_2)=\sum _{n\in \mathbb Z} (\operatorname P_n f_1)(\operatorname P_{n+k} f_2).
 \end{equation}
 It is straight forward to see that 
 \begin{equation}  \label{e.operatorDbound}
 \norm \operatorname D_k(f_1,f_2).p. {}\lesssim{}  \sup _{I\in \mathcal D } 
 \frac{ \abs { \ip f_1,h_I ,} } { \sqrt{\abs{I}} } \norm f_2.p. 
 \end{equation}
 In particular, these paraproducts admit a bound on their operator norms 
 that is strictly smaller than the $\operatorname {BMO}$ norm.

 \subsubsection{The Dyadic Reisz Potential and Commutator} 
 
 Consider a dyadic analog of the Reisz Potentials given by 
 \begin{equation}  \label{e.dyadicreisz}
 \operatorname I_\alpha f {}\eqdef{} \sum _{I\in \mathcal D} \frac{\ip f, \mathbf 1 _{I} , } {\abs I ^{\alpha}} \mathbf 1 _{I},
 \qquad 0<\alpha<1.
 \end{equation}
 This operator enjoys the same mapping properties of the continuous Reisz potentials, a fact we shall take for granted.

 And the continuous 
 versions can be recovered from the dyadic models by an appropriate averaging procedure. 
 This point of view is nicely illustrated in the article of Petermichl \cite{MR2000m:42016}, 
 in which the Hilbert transform is recovered from a dyadic model.

 We discuss the proof of a dyadic version of the Theorem of Chanillo, Theorem~\ref{t.chanillo}
 \begin{theorem} \label{t.dyadicchanillo}
 For $0<\alpha<1$,  $1-\alpha+\frac1q=\frac1p$ and $1<p<q<\infty$ we have 
 \begin{equation}  \label{e.dyadicchanillo}
 \norm [\operatorname M_b,\operatorname I_\alpha].p\to q .\simeq \norm b.\operatorname{BMO}. .
 \end{equation}
 \end{theorem}

 Indeed concerning the upper bound on the commutator, the main point  is this: The commutator 
 $[\operatorname I_\alpha,\operatorname M_b]$  is a linear combination of the four terms  
 \begin{gather}  \label{e.commutator=}
 		\operatorname B(b, \cdot)\circ \operatorname I_\alpha, \qquad 
		\operatorname I_\alpha\circ\operatorname D_0(b, \cdot)  ,
\\  \label{e.commutator=1}
 \operatorname D_0(b,\cdot  )\circ\operatorname I_\alpha,\qquad \sum _{k=1}^\infty 2 ^{-k(1-\alpha)} \operatorname
D_k(b,\cdot  )\circ\operatorname I_\alpha
 \end{gather}
 These operators are defined in (\ref{e.Haar-para-one}) and (\ref{e.operatorD}).
 Therefore, the upper bound on the commutator is an immediate consequence of those for the Reisz potentials, 
 and the corresponding paraproducts.

\medskip 

Observe that our Reisz potential, applied to a Haar function, has an explicit form. 
\begin{equation}  \label{e.reiszhaar}
\operatorname I_\alpha h_I=c_\alpha \abs I ^{1-\alpha} h_I , 
\end{equation}
for a choice of constant $c_\alpha=\sum _{n=1}^\infty 2 ^{-n(1-\alpha)}$.  In addition, for a dyadic interval $J$, we have 
\begin{equation} \label{e.reiszindicator}
\operatorname I_\alpha \mathbf 1 _{J}=\sum _{K\supset J} \frac {\abs J}{\abs K ^{\alpha} } \mathbf 1 _{K}.
\end{equation}
For later use, observe that $\operatorname I_\alpha \mathbf 1 _{J}$ equals $(1+c_\alpha) \abs{ J} ^{1-\alpha}$ 
on the interval $J$.

We can then compute the commutator, with multiplying function $h_I $ applied to another Haar function $h_J$
\begin{align} \nonumber 
[ \operatorname M _{h_I},\operatorname I_\alpha ]h_J&{}={} c _{\alpha} 
	h_I \abs J ^{1-\alpha} h_J-\operatorname I _{\alpha}(h_I\cdot h_J)
\\  \label{e.cases} 
&{}={}\begin{cases} 
0 & J\subsetneq  I
\\
{ c_\alpha}{\abs I} ^{-\alpha} \mathbf 1 _{I}- \abs I ^{-1}\operatorname I_\alpha( \mathbf 1 _{I})
		& I=J
\\
c_\alpha \, h_J(I)\{ \abs J ^{1-\alpha}-\abs I ^{1-\alpha} \} h_I & I\subsetneq J.
\end{cases}
\end{align}
And in the case that $I\subsetneq J$, 
note that $h_J$ takes exactly one value on $I$, which is denoted as $h_J(I)$.

 \medskip 
 
 We expand $[\operatorname M_b,\operatorname I_\alpha]f$ as a double sum over Haar functions.  The 
 leading term in the case of $I\subsetneq J$ in (\ref{e.cases}) gives us 
 \begin{align*}
 c _{\alpha}\sum _{I\in \mathcal D}\sum _{I\subsetneq J} {\ip b,h_I,} \ip f ,h_J, 
 \abs J ^{1-\alpha} h_J(I) h_I
 {}={}
 \operatorname B(b,\operatorname I_\alpha f)
 \end{align*}
 which is the first term in (\ref{e.commutator=}).
 
 For the second term in (\ref{e.cases}) in the case of $I\subsetneq J$,
 given a dyadic interval $I$ and integer $k>0$, let $I_k$ denote 
 the dyadic interval that contains $I$ and has length $\abs{I_k}=2^k\abs I $.  Note that 
 \begin{align*}
 \sum _{I\in \mathcal D} {\ip b,h_I,} \ip f ,h_{I_k}, 
 \abs I ^{1-\alpha} h_{I_k}(I) h_I=2 ^{-k(1-\alpha)}\operatorname D_k(b,\operatorname I_\alpha f). 
 \end{align*}
 This leads to the second half of (\ref{e.commutator=1}). 
 
 Consider the leading term in the case of $I=J$ in (\ref{e.cases}).  It gives us 
 \begin{equation*}
 \sum _{I\in \mathcal D} {\ip b,h_I,}  \ip f ,h_I, 
 \abs I ^{-\alpha} \mathbf 1 _{I}=\operatorname D_0(b,\operatorname I_\alpha f)
 \end{equation*}
 which is the first half of (\ref{e.commutator=1}).
 
 Consider the second term in the case of $I=J$ in (\ref{e.cases}).  It gives us 
 \begin{equation*}
 \sum _{I\in \mathcal D} {\ip b,h_I,}  \ip f ,h_I, 
 \abs I ^{-1} \operatorname I_\alpha \mathbf 1 _{I}=\operatorname I_\alpha\circ\operatorname D_0(b,f)
 \end{equation*}
 which is the second half of (\ref{e.commutator=}).  Our proof of  the upper bound on 
 the commutator norm in Theorem~\ref{t.dyadicchanillo} is finished. 
 
 \bigskip 
 
 Let us discuss the proof of the lower bound.  We can take $b\in \operatorname {BMO}$ of norm one. 
 Fix an interval $J$ so that 
 \begin{equation*}
 \sum _{I\subset J} \abs{\ip b, h_I,}^2\ge{} \tfrac12 \abs J. 
 \end{equation*}
 It is important to observe that by the John Nirenberg estimates we have 
 \begin{align*}
 1&{} {}\lesssim{} \abs{ J} ^{-1/p} \NOrm \sum _{I\subset J} \ip b, h_I,h_I.p.
 \\& {}\le{} \abs{ J} ^{-1/q} \NOrm \sum _{I\subset J} \ip b, h_I,h_I.q.
 \\& {}\lesssim{}1.
 \end{align*}

 We obtain a lower bound on the $L^q$ norm of the commutator applied to $\mathbf 1 _{J}$.  
 Write the function $b$ as $b=b'+b''$  where $b'=\sum _{\abs I\le\abs J} \ip b, h_I,h_I$. 
 \begin{align*}
 [ \operatorname M_b,\operatorname I_\alpha]\mathbf 1 _{J}& {}={} 
 	b(\operatorname I_\alpha \mathbf 1 _{J})-\operatorname I_\alpha(b \mathbf 1 _{J})
	\\&{}={} 
	b'\operatorname I_\alpha \mathbf 1 _{J}-\operatorname I_\alpha b'- b''{}\operatorname I_\alpha \mathbf 1 _{J}	
		{}+{} b''(J)\operatorname I_\alpha\mathbf 1 _{J}.
 \end{align*}
 Notice that $b''$ takes a single value on $J$, and that the last two terms cancel on that interval.  Thus, 
 \begin{equation} \label{e.lowercommutator} 
 \norm  [ \operatorname M_b,\operatorname I_\alpha]\mathbf 1 _{J}.q.\ge{} 
  \norm b'\operatorname I_\alpha \mathbf 1 _{J}-\operatorname I_\alpha b'.L^q(J).
 \end{equation}
  Taking the explicit formulas (\ref{e.reiszhaar})  and (\ref{e.reiszindicator}) into account, we see that the last term 
  above is at least a constant times 
  \begin{align*}
  \abs J ^{1-\alpha}\NOrm \sum _{I\subset J} \ip b, h_I,h_I .q.& {}\gtrsim{}
  	\abs J ^{1-\alpha+\tfrac1q-\tfrac1p}\NOrm \sum _{I\subset J} \ip b, h_I,h_I .p.
	\\ 
	& {}\gtrsim{} \abs J ^{1-\alpha+\tfrac1q }
	\\
	&{}={} \abs J ^{\tfrac1p}. 
  \end{align*}
  It follows that this commutator admits a universal lower bound on its $L^p\mapsto L^q$ norm, assuming 
  that the $\operatorname {BMO}$ of the function $b$ is one.  The proof is complete. 
  
  \smallskip 
  
  It is of interest to provide another proof of the lower bound.     Let us begin by establishing the lower bound 
  \begin{equation}  \label{e.firstlower}
  \norm [\operatorname M_b, \operatorname I_\alpha] .p\to q. 
  	{}\gtrsim{} \sup _{I\in \mathcal D} \frac{ \abs{\ip b,h_I,}} {\sqrt {\abs I}} . 
  \end{equation}
  Indeed, apply the commutator to the Haar function $h_I$, 
  \begin{align*}
  [\operatorname M_b, \operatorname I_\alpha]h_I& {}=c_\alpha b\cdot h_I-\operatorname I_\alpha (b h_I)
  \\
  &{}={}  -\ip b,h_I, \abs I ^{-1/2+1-\alpha}\mathbf 1 _{I}-c_\alpha
  	\sum _{J\subsetneq  I} h_I(J)\ip b,h_I, \frac{\abs I}{\abs J} 
  \end{align*}
  By the Littlewood Paley inequality for Haar functions, the latter term 
  can be ignored in providing a lower bound on the $L^q$ norm. We can then estimate 
  \begin{align*}
  \norm [\operatorname M_b, \operatorname I_\alpha]h_I .q. & {}\gtrsim{} 
  		\frac{\abs{ \ip b, h_I,} }{\sqrt{\abs I}} \abs I ^{1-\alpha} \norm \mathbf 1 _{I}.q.
	\\& {}\gtrsim{}
	 \frac{\abs{ \ip b, h_I,} }{\sqrt{\abs I}} \abs I ^{\frac1p}.  
  \end{align*}
  This proves (\ref{e.firstlower}).

  Now, in seeking to prove the lower bound, we can assume that $\norm b.\operatorname {BMO}.=1$, while 
  \begin{equation*}
  \sup _{I\in \mathcal D} \frac{\abs{ \ip b, h_I,} }{\sqrt{\abs I}}<\eta
  \end{equation*}
  where $\eta>0$ is a small absolute constant to be chosen.   
  
  Recall that the paraproducts  $\operatorname D_k$ have an upper bound 
  on their norm  given in  (\ref{e.operatorDbound}).   As well, we have shown that the 
  commutator $[ \operatorname M_b, \operatorname I_\alpha]$ as a sum of the terms in 
  (\ref{e.commutator=})---(\ref{e.commutator=1}).  Notice that for all of these terms, save one, 
  we have an upper bound on their norm of an absolute constant times $\eta$. 
 
  The one term that this does not apply to is $\operatorname B(b,\cdot)\circ \operatorname I_\alpha$. 
  But, it is very easy to see that 
  \begin{equation*}
  \norm \operatorname B(b,\cdot)\circ \operatorname I_\alpha.p\to q. {}\gtrsim{}c>0. 
  \end{equation*}
  Indeed, just apply the commutator to $\mathbf 1 _{J}$ for dyadic intervals $J$.
  And so, for $\eta>0$ sufficiently small, we see that 
  $\norm [\operatorname M_b, \operatorname I_\alpha].p\to q.>c/2$. 
 
 \section{Higher Parameter Commutators and Paraproducts } 
 
 We work in the setting of more variables, so that functions $f$ are defined on $\mathbb R^d$.
 Set $\operatorname I _{\alpha,j}$ to be the Reisz potential as defined in (\ref{e.dyadicreisz}), 
 applied in the $j$th coordinate.  For a sequence of choices of $0<\alpha_j<1$, observe that 
 the operator 
 \begin{equation*}
 \operatorname I_{\alpha_1}\circ \cdots \circ \operatorname I_{\alpha_d}
 \end{equation*}
 will map $L^p $ to $L^q$ provided $1-\sum _{j=1}^d \alpha_j+\frac1q=\frac1p$, and $1<p<q<\infty$.
 One uses the one parameter result in each coordinate seperately.

 Our main result is 
 
 \begin{theorem} \label{t.reiszmany}
  Let $0<\alpha_j<1$,  $1-\sum _{j=1}^d \alpha_j+\frac1q=\frac1p$, and $1<p<q<\infty$ we have 
  \begin{equation}  \label{e.reiszmany}
 \norm   [\cdots[\operatorname M_b,\operatorname I_{\alpha_1,1}],\cdots , \operatorname I_{\alpha_d,d}].p\to q. 
 {}\lesssim{} \norm b.\operatorname{BMO}_d.
  \end{equation}
 \end{theorem} 
 
The strategy of appealing to sharp function estimates has well known difficulties in 
the higher parameter setting,\footnote{R.~Fefferman \cite{MR90e:42030} has found a partial substitute 
for the sharp function in two parameters.}
and so we adopt the strategy  given in the previous section  in the one parameter setting. 
We recall the necessary 
results for the paraproducts in the higher parameter setting, and then detail the proof of 
the Theorem above.

Notice that this is not a full extension of Chanillo's result as we do no claim that the two  norms are comparable. 
We comment on this in more detail in Section~\ref{s.lower} below.

\subsection{Higher Parameter Paraproducts} 

Let $\mathcal R {}\eqdef{}\otimes _{j=1}^d \mathcal D $ denote the dyadic rectangles in $\mathbb R^d$. 
When needed, we will write such a rectangle as $R=\otimes _{j=1}^dR_j$. 
\begin{equation}  \label{e.haardef}
h_R(x_1,\ldots,x_d)=\prod_{j=1}^d h_{R_j}(x_j). 
\end{equation}
And, by $h_R^0$ we mean $h_R $.  The other distinguished 
function of this type is $ h^1_R=\abs{ h_R}$. 
 
The simplest higher parameter paraproduct, and the only one needed for this paper, is 
\begin{equation}  \label{e.dparaproducts}
\operatorname B (f_1,f_2) {}\eqdef{}\sum _{R\in \mathcal R} \frac{\ip f_1,h_R,}{\sqrt{\abs R}} 
\ip f_2,h_R^1, \, h_R . 
\end{equation}
 The principal fact about these paraproducts is this. 
 
 \begin{theorem} \label{t.dparaproducts}
 We have 
 \begin{equation}  \label{e.dparaproductsle}
 \norm \operatorname B (b,\cdot) .p. {}\simeq{}\norm b.\operatorname{BMO}_d.,\qquad 1<p<\infty. 
 \end{equation}
 \end{theorem} 
 
 In these inequalities,  the $\operatorname{BMO}_d $ space is the dual to product $H^1$, as identified by 
 S.-Y.~Chang and R.~Fefferman.  Specifically, 
 \begin{equation}  \label{e.BMOd}
 \norm b.\operatorname{BMO}_d.=\sup \Bigl[ \frac 1{\abs U}\sum _{R\subset U}\abs{\ip f,h_R,}^2\Bigr] ^{1/2}
 \end{equation}
 It is essential that in this definition, the supremum be formed over all open sets $U\subset \mathbb R^d$ of 
 finite measure.

  We caution the reader that the Theorem above does not include the full range of 
  multiparameters paraproducts.\footnote{The  presence of the full range of paraproducts 
 is the source of part of the difficulties in   Ferguson and Lacey 
 \cite{MR1961195} and Lacey and Terwelleger \cite{witherin}.}
 For more infomation about this theorem, see 
    Journ\'e \cite{MR836284}. More recently, see 
  Muscalu, Pipher, Tao and Thiele, \cites{camil,camil2} for certain extensions of the Theorem above. 
  Also see Lacey and Metcalfe \cite{laceymetcalfe}.

 \begin{proof} 
 The proof we will give will 
	rely upon the structure of the Hardy and $\operatorname {BMO}$ space, and the interpolation theory 
	for this pair of spaces. 
 
 It is efficient to establish appropriate end point estimates for the dual to this operator.  
 Fix $b\in \operatorname {BMO}_d$ of norm one.  We establish that the dual 
 operator $\operatorname B^\ast $ 
 maps $H^1\mapsto L^1$ and $L^\infty\mapsto \operatorname {BMO}_d$.  
 An interpolation argument will complete the proof. 
 
 For the $H^1$ estimate, we use the atomic theory, as given in \cite{MR86g:42038}.  
  Recall that an $H^1 $ atom is a function $\alpha$ with Haar support in a 
  set $A$ of finite measure, that is 
  \begin{equation} \label{e.atom1}
  \alpha=\sum _{R\subset  A} \ip \alpha,h_R, h_R .
  \end{equation}
  Moreover, it satisfies the size condition $\norm \alpha.2.\le{}\abs A ^{-1/2}$.  
  Every element $f\in H^1$ admits a representation $f=\sum _{j} c_j\alpha_j$ where each  $\alpha_j$ is an 
  atom, $c_j$ is a scalar, and 
  \begin{equation} \label{e.atom2}
  \norm f. H^1.\simeq\sum _j \abs{ c_j}.
  \end{equation}
 
 Observe that 
\begin{align*}
 \norm \operatorname B ^{\ast }(b,\alpha).1.&{}={} 
 	\sum_{R\subset A} \abs{ \ip b,h_R,\ip \alpha,h_R,}
\\&{}\le{} 
	\Bigl[  \sum_{R\subset A} \abs{ \ip b,h_R,}^2 \sum_{R\subset A} \abs{ \ip \alpha,h_R,}^2 
		\Bigr] ^{1/2}
\\&{}\le{}[\abs{ A}\abs{ A} ^{-1}] ^{1/2}=1. 
\end{align*}
 Thus, by (\ref{e.atom2}), it is clear that we have the $H^1\mapsto L^1$ estimate. 
 
 For the other estimate, fix a function $f\in L^\infty$ of norm one.  And take a set $U\subset \mathbb R^d $
 of finite measure.   Let us set 
 \begin{equation*}
 F_U {}\eqdef{}\sum _{R\subset U } \frac{\ip b,h_R,}{\sqrt{\abs{ R}} } \ip f,h_R, h^{1}_R.
 \end{equation*}
 Then, appealing to the definition of $\operatorname {BMO}$, it is the case that  
 \begin{align*}
 \norm F_U.1. & {}\le{} \sum_{R\subset U} \abs{ \ip b,h_R,\ip f,h_R,}
 \\ &{}\le{} 
 \Bigl[  \sum_{R\subset U} \abs{ \ip b,h_R,}^2 \sum_{R\subset U} \abs{ \ip f,h_R,}^2 
		\Bigr] ^{1/2}
\\&\le{} \abs{ U}.
 \end{align*}
 As this estimate is uniform over all choices of $U$, 
 it is a reflection of the John Nirenberg estimate in the multiparameter setting \cite{MR82a:32009} 
 that this implies that 
\begin{equation*}
 \norm F_U.2. {}\lesssim{}\abs U ^{1/2}. 
 \end{equation*}

  Let us observe that for  rectangles $R\subset U$ and $S\not\subset U$, we necessarily have 
 $\ip h_R,h^{1}_S,=0$.  Therefore, we have 
 \begin{align*}
 \sum_{R\subset U} \abs{ \ip {\operatorname B ^{\ast }(b,f)} , h_R, } ^2 
  	& {}={}\sum_{R\subset U} \abs{ \ip F_U,  h_R, } ^2 
\\&{}\le{} 
	\norm F_U .2.^2
\\& {}\lesssim{} 
	\abs U. 
 \end{align*}
 This proves the $L^\infty \mapsto \operatorname {BMO}_d$ bound.

 \end{proof}

  \subsection{A Secondary Result on Paraproducts} 
  
  For the proof of our main theorem, an estimate on certain paraproducts is needed.  These paraproducts 
  are of a secondary nature.  In particular, we can give an upper bound on their norm that is 
  strictly smaller, in general, than the $\operatorname{BMO}_d$ norm.  
  
  The paraproducts are easiest to define in terms of tensor products of operators.
  Let $\mathbb B$ and $\mathbb D$ be a  partition of the coordinates 
  $\{1,\cdots,d\}$, and define 
  \begin{equation}  \label{e.tensor2}
  \operatorname E _{\mathbb B}(b,\cdot) {}\eqdef{}\mathop\otimes _{j\in \mathbb B}\operatorname B_j(b,\cdot)\mathop\otimes _{j\in \mathbb D} 
  \operatorname D_{v(j),j}(b,\cdot)
  \end{equation}
  Here, $v(j)$ is a non negative integer, and $\operatorname D _{v(j),j}$ is the operator $\operatorname D _{v(j)}$ 
  as defined in (\ref{e.operatorD}) acting in the $j$th coordinate.

  \begin{proposition}\label{p.paraproductless}
  The paraproducts $\operatorname E _{\mathbb B}$ admit the bound 
  \begin{equation}  \label{e.paraproductless}
  \norm \operatorname E _{\mathbb B}(b,\cdot) .p\to p. {}\lesssim{}
  \begin{cases} 
   \norm b . \operatorname{BMO} _{d,\mathbb  B} .   & 1<p\le{} 2 
  \\ 
  \norm b . \operatorname{BMO} _{d,\mathbb  B} . ^{2/p} \norm b . \operatorname{BMO} _{d} . ^{1-2/p}   &  2<p< \infty. 
  \end{cases} 
  \end{equation}
  \end{proposition}
  
  In this Proposition, the  norm $\norm \cdot.\operatorname{BMO} _{d,\mathbb  B}. $ is defined in terms of collections of rectangles $\mathcal S$ which 
  are restricted in the following way.  Say that $\mathcal S$ is of type $\mathbb B$ if the rectangles in 
  $\mathcal S$ have a union with finite measure, and for all $R,R'\in \mathcal S$,
  and coordinates $j\not\in \mathbb B $ 
  we have $R_j=R'_j$.   Thus, only the coordinates in $\mathbb B$ are permitted to vary.  Then define 
  \begin{equation}  \label{e.BMOB}
  \norm b.\operatorname{BMO} _{d, \mathbb B}. {}\eqdef{}\sup _{\mathcal S} \Bigl[\Abs{\bigcup _{R\in \mathcal S} R } ^{-1}
  \sum _{R\in \mathcal S} \abs{\ip f,h_R,}^2  \Bigr] ^{1/2}
  \end{equation}
  where the supremum is over all collections of rectangles of type $\mathbb B$. 
 Clearly, this norm is strictly smaller than that of $\operatorname{BMO}_d $.   And an example of Carleson \cite{carleson-example}, 
 and published in \cite{MR90e:42030}, shows that these norms are essentially smaller than the $\operatorname {BMO }$ norm. 
 (The use of norms of these types are illustrated in \cite{witherin} and \cite{journesurvey}.) 
 
 \begin{proof}
  We proceed to the proof. The case of the cardinality of $\mathbb B$ is 
  full, that is  equal to $d$, is contained in Theorem~\ref{t.dparaproducts}.  
  
  Now assume that the cardinality of $\mathbb B$ is not full. The $L^2$ case of (\ref{e.paraproductless})
  follows immediately, as we are forming the tensor product of operators $\operatorname D_v$ that act 
  on a family of orthogonal spaces.   
  
  We should take care to consider the form of the operator $\operatorname E _{\mathbb B}(b,f) $.  
  For a dyadic rectangle $R$, and coordinate $j $, let $\widetilde R_j=R_j$ if $j\not\in \mathbb D$, and 
  otherwise take this to be the dyadic 
  interval that contains $R_j$ and has length $2 ^{v(j)}\abs {R_j} $.   Let $\widetilde R=\otimes \widetilde R_j$. 
  Then, the operator in question is 
  \begin{equation*}
  \operatorname E _{\mathbb B}(b,f)=\sum _{R} \varepsilon_R 
  	\frac{ \ip b,h_R,}{\sqrt{\abs{\widetilde R}}} \ip f, h ^{\epsilon}_{\widetilde R } ,h ^{\epsilon}_{\widetilde R } 
  \end{equation*}
  Here, $\epsilon\in\{0,1\}^d$ is equal to $1$ for those coordinates in 
  $\mathbb B$, and is zero otherwise, and 
  \begin{equation*}
h ^{\epsilon}_{\widetilde R }=\prod _{j=1}^d h ^{\epsilon_j} _{\widetilde R_j}.
\end{equation*}
  The coefficient $\varepsilon_R$ is a choice of sign. (To be specific, the value of $\varepsilon_R$ 
  is the product of the signs $\operatorname {sgn}(h _{\widetilde R_j}(R_j)$ over those $j\in \mathbb D$ 
  such that $v(j)>0$.)

  It is natural to exploit the availible $L^2$ estimate by establishing the boundedness of $\operatorname E _{\mathbb B}$
  as a map from   $H^1 \mapsto H^1$. Recall the definition of an $H^1$ atom in (\ref{e.atom1}) and (\ref{e.atom2})

  Since $R\subset \widetilde R$, it 
  is clear that $\operatorname E _{\mathbb B}$  applied to an atom has the same Haar support.  And by the $L^2$ 
  bound, 
  \begin{equation*}
  \norm \operatorname E _{\mathbb B} \alpha.2. {}\lesssim{}\norm b.\operatorname{BMO} _{d,\mathbb B}. \abs A ^{-1/2}
  \end{equation*}
  This proves the bound at $H^1$. And by interpolation, we deduce the result for $1<p<2$. 
  
  At the other endpoint, we prove that 
  \begin{equation*}
  \norm \operatorname E _{\mathbb B}(b,\cdot). L^\infty\mapsto \operatorname{BMO} _{d}. {}\lesssim{} \norm b .\operatorname{BMO}_d. 
  \end{equation*}
  Indeed, take $f\in L^\infty$ of norm one, and a set $U\subset \mathbb R^d$ of finite measure. 
  Then, 
  \begin{align*}
  \sum _{\widetilde R\subset U} \abs{ \ip b, h_R,}^2 \frac{\abs{ \ip f, h ^{\varepsilon}_{\widetilde R},}^2}{\abs {\widetilde R}} 
  {}\le{} \norm f.\infty.^2\sum _{\widetilde R\subset U} \abs{ \ip b, h_R,}^2\le{} \norm f.\infty.^2 \norm b.\operatorname{BMO}_d.^2
  \end{align*}
  This proves the inequality at $L^\infty$, and interpolation will prove the bound for $2<p<\infty$. 
  
  \end{proof}

 \subsection{The Proof of Theorem~\ref{t.reiszmany}.} 
 
 The same proof strategy as in one dimension is used.  We expand the commutator  as 
 a double sum over Haar functions.  In so doing, we use (\ref{e.cases}).  Observe that 
 for two rectangles $R$ and $S$, we have 
 \begin{align*}
  [\cdots[\operatorname M_{h_R},\operatorname I_{\alpha,1}],\cdots , \operatorname I_{\alpha,d}] h_S
  &{}={} \prod_{j=1}^d [\operatorname M _{h_{R_j}},\operatorname I_{\alpha_j}]h _{S_j}
  \\&{}={}  0
 \end{align*}
 if for any coordinate $j$ we have $S_j\subsetneq R_j$.   Assuming that this is not the case, 
 we see that one of two terms can arise in each coordinate, depending upon 
 $S_j=R_j$ or $R_j\subsetneq  S_j$, as described in (\ref{e.cases}).  In this way, we 
 expand the commutator as sum of paraproduct operators.

  These operators  are as in (\ref{e.tensor2}):
 \begin{equation} \label{e.allparas} 
 2 ^{-v(1-\alpha)}\operatorname E _{\mathbb B}(b,\cdot) \circ \otimes _{j=1}^d\operatorname I _{\alpha_j,j},\qquad 
  2 ^{-v(1-\alpha)} \otimes _{j=1}^d\operatorname I _{\alpha_j,j} \circ \operatorname E _{\mathbb B}(b,\cdot).
 \end{equation} 
 We permit the subset $\mathbb B\subset \{1,\ldots,d\}$ to 
 vary over all possible subsets.   Associated to the complementary set $\mathbb D=\{1,\ldots,d\}-\mathbb B$ 
 is a vector $v=\{v(j)\}\in \mathbb N^{\mathbb D}$, and we set $v=\sum _{j\in \mathbb D}v(j)$. 
 
 According to Proposition~\ref{p.paraproductless} and the obvious bound on the Reisz potential,
 each of these terms has $L^p\mapsto L^q$ norm of at most $2 ^{-v(1-\alpha)} \norm b.\operatorname{BMO}_d.$. 
 And these estimates are summable over all choices of $\mathbb B$, $\mathbb D$, and choices of 
 integers $\{v(j)\mid j\in \mathbb D\}$. This completes the proof of the upper bound.

 \subsection{Concerning Lower Bounds on the Commutator Norm.}  \label{s.lower} 
 In the one parameter case, to provide a lower bound on the norm of the paraproduct 
 $\operatorname B(b,\cdot)$, as defined in (\ref{e.Haar-para-one}), it suffices to test it against an 
 indicator of a dyadic interval.  Moreover, one trivially has 
 \begin{equation*}
 \norm \operatorname  B(b,\operatorname I _{\alpha}\mathbf 1 _{J}).q. {}\gtrsim{} c_\alpha \abs J ^{1-\alpha} 
 \norm \operatorname  B(b,\mathbf 1 _{J}).q..
 \end{equation*}

 In higher parameters, the situation is far less obvious.  We can establish 
 
 \begin{proposition}\label{p.bmorec} We have the inequality 
 \begin{equation*} 
 \norm  [\cdots[\operatorname M_b,\operatorname I_{\alpha,1}],\cdots , \operatorname I_{\alpha,d}].p\to q.
 {}\gtrsim{}\sup _{S } \Bigl[ \abs{ S} ^{-1} \sum _{R\subset S} \abs {\ip b,h_R,}^2\Bigr] ^{1/2}.
 \end{equation*}
 Here, the supremum is formed over all  dyadic rectangles $S$. 
  \end{proposition}

  We omit the proof, which depends upon an iteration of the argument that lead to (\ref{e.lowercommutator}).  
  
 This norm sometimes referred to as ``rectangular $\operatorname {BMO}$,'' denoted as $\norm \cdot.\operatorname 
 {BMO(rec)}.$.  
 But as is well known, this norm is essentially smaller than the $\operatorname {BMO}$ norm, and so 
 this proposition is not enough to prove the complete analog of Chanillo's theorem.  
 
 Continuing this line of thought, in two dimensions (and only two dimensions), it is easy to see that we have 
 \begin{equation*}
 \norm b .\operatorname {BMO(rec)}.={} \sup _{\mathbb B=\{1\},\{2\}} \norm b.\operatorname {BMO} _{2,\mathbb B}.
 \end{equation*}
 From this, our expansion of the commutator, and Proposition~\ref{p.paraproductless}, we see that we have the 
 estimate 
 \begin{align*}
 \norm [[\operatorname M_b,\operatorname I _{\alpha_1,1}],\operatorname I _{\alpha_2,2}]&. p\to q.
 {}\gtrsim{}
\\&\max( \norm b.\operatorname {BMO(rec)}. ^{\varepsilon(\alpha,p)}, \ 
 \norm \operatorname B(b,\cdot)\circ \operatorname I _{\alpha_1,1}
\otimes  
 \operatorname  I _{\alpha_2,2}+\operatorname I _{\alpha_1,1}
\otimes  
 \operatorname  I _{\alpha_2,2}\circ \operatorname B(b,\cdot).p \to q.  )
 \end{align*}
 Here, $\varepsilon(\alpha,p) $ is a positive exponent.  
 It is natural to then assume that the rectangular norm is small, and argue that the 
 other norm is big, but there is a problem 
 with continuing this line of thought.

 To provide a lower bound on the norms of 
 the paraproducts $\operatorname B (b,\cdot)$ as defined in (\ref{e.dparaproducts}), 
 we need to apply this operator to $\mathbf 1 _{U}$, for arbitrary sets $U\subset \mathbb R^d$ of 
 finite measure.   We would then like for an inequality of this form to be true:
 \begin{equation*}
 \norm \operatorname I_{\alpha_1,1}\circ\cdots\circ \operatorname I_{\alpha_d,d} \mathbf 1 _{U}.q. {}\gtrsim{} 
 \abs U ^{1/p}, \qquad U\subset \mathbb R^d.
 \end{equation*}
  This holds for $U$  a rectangle, but not in general. Indeed, consider $U=\bigcup_{n=1}^n R_n$, where $R_n$ are rectangles that 
 are of the same dimension, but very widely seperated.  So seperated that we can estimate 
 \begin{align*}
 \norm \operatorname I_{\alpha_1,1}\circ\cdots\circ \operatorname I_{\alpha_d,d} \mathbf 1 _{U}.q.
 & {}\simeq \Bigl[ \sum_{n=1}^N 
 \norm \operatorname I_{\alpha_1,1}\circ\cdots\circ \operatorname I_{\alpha_d,d} \mathbf 1 _{R_n}.q.^q 
 \Biggr] ^{1/q}
 \\ 
 & {}\simeq{} N ^{1/q}\abs{ R} ^{1/p}
 \\
 &{}\simeq{} N ^{1/q-1/p}\abs U ^{1/p}.
 \end{align*}
 Since $1<p<q<\infty$, this last term is substantially smaller than $\abs{ U} ^{1/p}$ for $N$ large.


\smallskip 
 
 
 

 
 
  \begin{bibdiv} 
 \begin{biblist} 
 
 
 
 
 \bib{MR1934198}{article}{
    author={Auscher, P.},
    author={Hofmann, S.},
    author={Muscalu, C.},
    author={Tao, T.},
    author={Thiele, C.},
     title={Carleson measures, trees, extrapolation, and $T(b)$ theorems},
   journal={Publ. Mat.},
    volume={46},
      date={2002},
    number={2},
     pages={257\ndash 325},
      issn={0214-1493},
    review={MR1934198 (2003f:42019)},
}

\bib{journesurvey}{article}{
 author={Cabrelli, Carlos},
 author={Molter, Ursula},
 author={Lacey, Michael T.},
 author={Pipher, Jill},
 title={Variations of a Theme of Journ\'e},
 journal={To appear in Houston J. Math.,\href{http://www.arxiv.org/abs/math.CA/0310367}{arxiv:math.CA/0310367}},
 }   
 
 
 
 
 \bib{carleson-example}{article}{
	author={Carleson, L.},
	title={A counterexample for measures bounded on $H^p$ spaces for the bidisk},
	journal={Mittag-Leffler Rep. No. 7, Inst. Mittag-Leffler},
	year={1974},
	}


\bib{MR86g:42038}{article}{
    author={Chang, Sun-Yung A.},
    author={Fefferman, Robert},
     title={Some recent developments in Fourier analysis and $H\sp p$-theory
            on product domains},
   journal={Bull. Amer. Math. Soc. (N.S.)},
    volume={12},
      date={1985},
    number={1},
     pages={1\ndash 43},
      issn={0273-0979},
    review={MR 86g:42038},
}

\bib{MR82a:32009}{article}{
    author={Chang, Sun-Yung A.},
    author={Fefferman, Robert},
     title={A continuous version of duality of $H\sp{1}$ with {BMO} on the
            bidisc},
   journal={Ann. of Math. (2)},
    volume={112},
      date={1980},
    number={1},
     pages={179\ndash 201},
      issn={0003-486X},
    review={MR 82a:32009},
}


 \bib{MR84j:42027}{article}{
    author={Chanillo, S.},
     title={A note on commutators},
   journal={Indiana Univ. Math. J.},
    volume={31},
      date={1982},
    number={1},
     pages={7\ndash 16},
      issn={0022-2518},
    review={MR 84j:42027},
}


\bib{MR0380244}{article}{
    author={Coifman, R. R.},
    author={Meyer, Yves},
     title={On commutators of singular integrals and bilinear singular
            integrals},
   journal={Trans. Amer. Math. Soc.},
    volume={212},
      date={1975},
     pages={315\ndash 331},
    review={MR0380244 (52 \#1144)},
}

\bib{MR511821}{article}{
    author={Coifman, R. R.},
    author={Meyer, Yves},
     title={Commutateurs d'int\'egrales singuli\`eres et op\'erateurs
            multilin\'eaires},
  language={French, with English summary},
   journal={Ann. Inst. Fourier (Grenoble)},
    volume={28},
      date={1978},
    number={3},
     pages={xi, 177\ndash 202},
      issn={0373-0956},
    review={MR511821 (80a:47076)},
}

\bib{MR518170}{book}{
    author={Coifman, R. R.},
    author={Meyer, Yves},
     title={Au del\`a des op\'erateurs pseudo-diff\'erentiels},
  language={French},
    series={Ast\'erisque},
    volume={57},
 publisher={Soci\'et\'e Math\'ematique de France},
     place={Paris},
      date={1978},
     pages={i+185},
    review={MR518170 (81b:47061)},
}





\bib{MR54:843}{article}{
    author={Coifman, R. R.},
    author={Rochberg, R.},
    author={Weiss, Guido},
     title={Factorization theorems for Hardy spaces in several variables},
   journal={Ann. of Math. (2)},
    volume={103},
      date={1976},
    number={3},
     pages={611\ndash 635},
    review={MR 54 \#843},
}


\bib{MR2004a:42010}{article}{
    author={Cruz-Uribe, D.},
    author={Fiorenza, A.},
     title={Endpoint estimates and weighted norm inequalities for
            commutators of fractional integrals},
   journal={Publ. Mat.},
    volume={47},
      date={2003},
    number={1},
     pages={103\ndash 131},
      issn={0214-1493},
    review={MR 2004a:42010},
}




\bib{MR763911}{article}{
    author={David, Guy},
    author={Journ{\'e}, Jean-Lin},
     title={A boundedness criterion for generalized Calder\'on-Zygmund
            operators},
   journal={Ann. of Math. (2)},
    volume={120},
      date={1984},
    number={2},
     pages={371\ndash 397},
      issn={0003-486X},
    review={MR763911 (85k:42041)},
}


\bib{MR2002g:42016}{article}{
    author={Ding, Yong},
    author={Lu, Shanzhen},
    author={Zhang, Pu},
     title={Weak estimates for commutators of fractional integral operators},
   journal={Sci. China Ser. A},
    volume={44},
      date={2001},
    number={7},
     pages={877\ndash 888},
      issn={1006-9283},
    review={MR 2002g:42016},
}


\bib{MR2084076}{article}{
    author={Duong, Xuan Thinh},
    author={Yan, Li Xin},
     title={On commutators of fractional integrals},
   journal={Proc. Amer. Math. Soc.},
    volume={132},
      date={2004},
    number={12},
     pages={3549\ndash 3557 (electronic)},
      issn={0002-9939},
    review={MR2084076 (2005e:42046)},
}

\bib{MR90e:42030}{article}{
    author={Fefferman, Robert},
     title={Harmonic analysis on product spaces},
   journal={Ann. of Math. (2)},
    volume={126},
      date={1987},
    number={1},
     pages={109\ndash 130},
      issn={0003-486X},
    review={MR 90e:42030},
} 


\bib{MR86f:32004}{article}{
    author={Fefferman, R.},
     title={A note on Carleson measures in product spaces},
   journal={Proc. Amer. Math. Soc.},
    volume={93},
      date={1985},
    number={3},
     pages={509\ndash 511},
      issn={0002-9939},
    review={MR 86f:32004},
}

\bib{MR81c:32016}{article}{
    author={Fefferman, R.},
     title={Bounded mean oscillation on the polydisk},
   journal={Ann. of Math. (2)},
    volume={110},
      date={1979},
    number={2},
     pages={395\ndash 406},
      issn={0003-486X},
    review={MR 81c:32016},
}


 \bib{MR1961195}{article}{
    author={Ferguson, Sarah H.},
    author={Lacey, Michael T.},
     title={\href{http://www.arxiv.org/abs/math.CA/0104144}{A characterization of product {BMO} by commutators}},
   journal={Acta Math.},
    volume={189},
      date={2002},
    number={2},
     pages={143\ndash 160},
      issn={0001-5962},
    review={1 961 195},
}
 \bib{MR80j:42034}{article}{
    author={Janson, Svante},
     title={Mean oscillation and commutators of singular integral operators},
   journal={Ark. Mat.},
    volume={16},
      date={1978},
    number={2},
     pages={263\ndash 270},
      issn={0004-2080},
    review={MR 80j:42034},
}





\bib{MR836284}{article}{
    author={Journ{\'e}, Jean-Lin},
     title={Calder\'on-Zygmund operators on product spaces},
   journal={Rev. Mat. Iberoamericana},
    volume={1},
      date={1985},
    number={3},
     pages={55\ndash 91},
      issn={0213-2230},
    review={MR836284 (88d:42028)},
}



\bib{laceymetcalfe}{article}{
	author={Lacey, Michael T.},
	author={Metcalfe, Jason},
	title={Paraproducts in One and Several Parameters  },
	date={2005},
}

 

\bib{witherin}{article}{
	author={Lacey, Michael T.},
	author={Terwilleger, Erin},
	title={Little Hankel Operators and Product {BMO}},
	date={2004},
}






\bib{MR1783613}{article}{
    author={Lacey, Michael},
    author={Thiele, Christoph},
     title={A proof of boundedness of the Carleson operator},
   journal={Math. Res. Lett.},
    volume={7},
      date={2000},
    number={4},
     pages={361\ndash 370},
      issn={1073-2780},
    review={MR1783613 (2001m:42009)},
}





\bib{nehari}{article}{
    author={Nehari, Z.},
     title={On bounded bilinear forms},
   journal={Ann. of Math. (2)},
    volume={65},
      date={1957},
    number={},
     pages={153\ndash 162},
      issn={},
    review={},
}


\bib{camil}{article}{
 author={Mucalu, Camil},
 author={Pipher, Jill},
 author={Tao, Terrance},
 author={Thiele, Christoph},
 title={Bi-parameter paraproducts},
 journal={\href{http://www.arxiv.org/abs/math.CA/0310367}{arxiv:math.CA/0310367}},
 }
 

 \bib{camil2}{article}{
 author={Mucalu, Camil},
 author={Pipher, Jill},
 author={Tao, Terrance},
 author={Thiele, Christoph},
 title={Multi-parameter paraproducts},
 journal={\href{http://www.arxiv.org/abs/math.CA/0411607}{arxiv:math.CA/0411607}},
 }
 
 







   

\bib{MR96j:42006}{article}{
    author={Paluszy{\'n}ski, M.},
     title={Characterization of the Besov spaces via the commutator operator
            of Coifman, Rochberg and Weiss},
   journal={Indiana Univ. Math. J.},
    volume={44},
      date={1995},
    number={1},
     pages={1\ndash 17},
      issn={0022-2518},
    review={MR 96j:42006},
}




\bib{MR1964822}{article}{
    author={Petermichl, S.},
    author={Treil, S.},
    author={Volberg, A.},
     title={Why the Riesz transforms are averages of the dyadic shifts?},
   journal={Publ. Mat.},
      date={2002},
    number={Vol. Extra},
     pages={209\ndash 228},
      issn={0214-1493},
    review={MR1964822 (2003m:42028)},
}



\bib{MR2000m:42016}{article}{
    author={Petermichl, Stefanie},
     title={Dyadic shifts and a logarithmic estimate for Hankel operators
            with matrix symbol},
  language={English, with English and French summaries},
   journal={C. R. Acad. Sci. Paris S\'er. I Math.},
    volume={330},
      date={2000},
    number={6},
     pages={455\ndash 460},
      issn={0764-4442},
    review={MR 2000m:42016},
}




\bib{MR1232192}{book}{
    author={Stein, Elias M.},
     title={Harmonic analysis: real-variable methods, orthogonality, and
            oscillatory integrals},
    series={Princeton Mathematical Series},
    volume={43},
 publisher={Princeton University Press},
     place={Princeton, NJ},
      date={1993},
     pages={xiv+695},
      isbn={0-691-03216-5},
    review={MR1232192 (95c:42002)},
}




 

 \bib{MR96k:47089}{article}{
    author={Youssfi, Abdellah},
     title={Regularity properties of commutators and ${\rm
            \operatorname{BMO}}$-Triebel-Lizorkin spaces},
  language={English, with English and French summaries},
   journal={Ann. Inst. Fourier (Grenoble)},
    volume={45},
      date={1995},
    number={3},
     pages={795\ndash 807},
      issn={0373-0956},
    review={MR 96k:47089},
}




 
 \end{biblist}    
 \end{bibdiv}
\end{document}